%
%
%
\input amstex.tex
\input amsppt.sty
\input pictex
\documentstyle{amsppt}
\NoBlackBoxes
\let\e\varepsilon
\let\phi\varphi
\define\N{{\Bbb N}}

\define\OO{{\Cal O}}

\def\ss{\subset}
\def\vv{\vert}

\let\i\infty

\def\g{\gamma}

\def\a{\alpha}
\def\b{\beta}
\def\q{\qed}
\def\f{\frac}
\def\qq{\qquad}

\def\int{\hbox{int}}

\topmatter
\title   A two--parameter control for contractive-like multivalued mappings  \endtitle

\author       Du\v{s}an Repov\v{s} \endauthor
\leftheadtext{Du\v{s}an Repov\v{s} }
\rightheadtext{A two--parameter control}

\address
Faculty of Mathematics and Physics, and
Faculty of Education, 
University of Ljubljana,
P.O.B. 2964, Ljubljana, Slovenia 1001
\endaddress
\email dusan.repovs\@guest.arnes.si
\endemail


\abstract
We propose a general approach to defining a contractive-like multivalued mappings $F$ which avoids any use of the Hausdorff distance between the sets $F(x)$ and $F(y)$. Various fixed point theorems are proved under a two-parameter control of the distance function $d_{F}(x)=dist(x,F(x))$ between a point $x \in X$ and the value $F(x) \ss X$. Here, both parameters are numerical functions. The first one $\a\,:[0,+\i)\rightarrow [1,+\i)$ controls the distance between $x$ and some appropriate point $y \in F(x)$ in comparison with $d_{F}(x)$, whereas the second one $\b\,:[0,+\i)\rightarrow [0,1)$ estimates $d_{F}(y)$ with respect to $d(x,y)$. It appears that the well harmonized relations between $\a$ and $\b$ are sufficient for the existence of fixed points of $F$. Our results  generalize several known fixed-point theorems. 
\endabstract

\endtopmatter

\document

\head {\bf 0. Introduction}\endhead

Fixed points theories for singlevalued mappings and  multivalued mappings of
metric spaces into themself are clearly closely related. As a rule, almost any fixed point
theorem for multivalued mappings goes back to
some fixed point theorem for singlevalued mappings.

Such  a correlation basically deals with the substitution of a given metric, say $d$,
on a space $X$ by the corresponding Hausdorff "metric" $H_{d}$ on the set of all closed subsets of $X$.
To illustrate this point, if one proves a fixed point theorem for a singlevalued mapping
$f:X\rightarrow X$ under some contractivity-type restriction as e.g.\,Browder's condition \cite{1,5,6}:
$$
d(f(x),f(y)) \leq \varphi(d(x,y))< d(x,y) \qq \qq \qq (*)
$$
then one can certainly  try to verify the existence of fixed points for a multivalued mapping
$F$ from $X$ into itself under the analogous assumption:
$$
H_{d}(F(x),F(y)) \leq \varphi(d(x,y))< d(x,y).        \qq \qq \qq (**)
$$

Recall that the inequality $H_{d}(A,B) < \e$ implies that each of the sets $A$ and $B$ is a subset of an open $\e-$neighborhood of the other set. The key goal of the present paper is to show that the proximity  of $F(x)$ and $F(y)$ with respect to $H_{d}$ is too restrictive for a successful construction of the Picard sequence of approximations $x_{n} \rightarrow x_{*}$ which would converge to a fixed point $x_{*}$ of $F$,\,\, $x_{*} \in F(x_{*})$.

Roughly speaking, there is no need to require that the {\it entire } set $F(x_{n})$ lies in an
$\e-$neighborhood $\OO_{\e}(F(x_{n+1}))$ of the set $F(x_{n+1})$ and, symmetrically that $F(x_{n+1}) \ss \OO_{\e}(F(x_{n}))$. It suffices to find for a chosen  $x_{n} \in X$, a point $x_{n+1} \in F(x_{n})$ such that the distance $d_{n}=d(x_{n},x_{n+1})$ is  "almost" equal to   $dist(x_{n}, F(x_{n}))$  and additionally, the distance $dist(x_{n+1}, F(x_{n+1}))$ is "less" than $d_{n}$. For the control of "nearness" of $d_{n}$ to $dist(x_{n}, F(x_{n}))$ we propose a numerical function
$\a:[0,+\i)\rightarrow [1,+\i)$, whereas the control of $dist(x_{n+1}, F(x_{n+1}))$ with respect to $d_{n}$ will be provided by a numerical function $\b:[0,+\i)\rightarrow [0,1)$.
Certain matching behaviour of control functions $\a$ and $\b$ guarantees the convergence of the sequence $\{x_{n}\}$. A standard verification shows that the limit of the sequence $\{x_{n}\}$ is a fixed point of the multivalued mapping $F$.

\head {\bf 1. Preliminaries}\endhead

For a metric space $(X,d)$, a point $x \in X$ and  a subset $A \ss X$ we denote by $\OO_{\e}(x)$ the open $\e-$neighborhood of $x$ and  by $\OO_{\e}(A)$ the open $\e-$neighborhood of $A$, i.e.\,\,$\OO_{\e}(A) = \bigcup \{ \OO_{\e}(x) \vv\,\,x \in A \}$. The distance between $x$ and $A$  is defined by $dist(x,A)= \inf \{d(x,y) \vv y \in A \}$, as usual. For a fixed closed-valued mapping $F$ of a metric space $(X,d)$ into itself we denote by $d_{F}(\cdot)$ the distance function which is defined by the equality
$
d_{F}(x)=dist(x,F(x)),\,\,\,x \in X.
$

For nonempty subsets $A \ss X$ and $B \ss X$ the {\it Hausdorff distance} between $A$ and $B$ is defined by $H_{d}(A,B)=\inf \{\e>0 \vv A\ss \OO_{\e}(B),\,\, B\ss \OO_{\e}(A)\}$, or equivalently, by $H_{d}(A,B)=\max \{ \sup \{dist(x,B)\vv x \in A \}, \sup \{dist(y,A)\vv y \in B \} \}$. In order to get the property $(H_{d}(A,B)=0 \,\,\Leftrightarrow \,\, A=B)$ one needs to use the Hausdorff distance only for closed subsets of $X$. It is well known that $H_{d}(\cdot,\cdot)$  is indeed a metric on the family $CB(X)$ of all closed bounded nonempty subsets of $X$. Moreover, for a complete metric space $(X,d)$ the {\it compact exponent} of $X$, i.e. the family of all nonempty subcompacta of $X$,\, is a complete metric space with respect to $H_{d}(\cdot,\cdot)$.

Let $k:[0,\i) \to [0;1)$ be any numerical function. A multivalued mapping $F$ which
associates to each point $x$  of a metric space $(X,d)$
some nonempty closed subset $F(x) \ss X$, is called  a {\it Hausdorff $k$-contraction} if
$ H_{d}(F(x),F(y)) \le  k(d(x,y))\cdot d(x,y),\,\,x,y \in X.$
Among others we extract two rather typical restrictions for the contractivity coefficient.
A numerical function $k:[0,\i) \to [0;1)$ is said to have the {\it Reich property} $(R)$ 
if \,\,$\lim\sup_{s \to t+0} k(s) \,\,<\,\,1$ for every $t > 0$  (cf. \cite{9, 10}),
and to have the {\it Mizoguchi-Takahashi property} $(MT)$ if the same inequality holds for every $t \geq 0$ (cf. \cite{8}).

The following is the key notion of the present paper.

\proclaim{Definition 1.1}  For numerical functions
$\a:[0,+\i)\rightarrow [1,+\i)$ and $\b:[0,+\i)\rightarrow [0,1)$, a closed-valued mapping $F$ of a metric space $(X,d)$ into itself is said to be an $(\a, \b)-$mapping if for each $x \in X$ there exists $y \in F(x)$ such that:
\itemitem{($A$)} $d(x,y) \leq \a(d(x,y))\cdot d_{F}(x)$; and
\itemitem{($B$)} $d_{F}(y) \leq \b(d(x,y))\cdot d(x,y)$.

If, in addition, $\a(t)\cdot \b(t) < 1,\,\, t > 0$,\, then $F$ is said to be an $(\a, \b)-$contraction.
\endproclaim

It easy to see that every Hausdorff $k$-contraction $F$ is an $(\a, k)-$contraction for every numerical function $\a:[0,+\i)\rightarrow (1,+\i)$ with $\a(t)\cdot k(t) < 1,\, t > 0$. In fact,
the assumption $(B)$ is true for an arbitrary  $y \in F(x)$ because $d_{F}(y) \leq H_{d}(F(x),F(y)) \leq k(d(x,y))\cdot d(x,y)=\b(d(x,y))\cdot d(x,y)$, whereas the assumption $(A)$ evidently holds for $d_{F}(x)=0 \Leftrightarrow x \in F(x) \Leftrightarrow y=x \Leftrightarrow d(x,y)=0$ and for $d_{F}(x)>0$ it is true for a suitable $y \in F(x)$ because $d_{F}(x) < \a(d(x,y))\cdot d_{F}(x)$. Hence, each fixed point theorem for a Hausdorff contraction can be considered as a special case of some fixed point theorem for an $(\a, \b)$-contraction.

We shall also need some technical notations.

\proclaim{Definition 1.2}
\itemitem{(1)} A numerical function $h:[0,\i) \to [0;1)$ is said to be essentially positive
if \break $\inf\{h(s)\,|\,\,s\geq a\}>0$ for every $a>0$;
\itemitem{(2)} A nonnegative (and non-identically
zero) function $h:X \to [0;1)$ is said to be
stably positive if the inequality $h(x)>0$ implies that $\inf\{h(y)\,|\, \,y\in \OO(x)\}>0$ for some neighborhood $\OO(x)$ of the point $x \in X$.
\endproclaim

To verify the implication $(x_{n} \to x_{*})\wedge(x_{n+1} \in F(x_{n}), n \in \N) \Longrightarrow x_{*} \in F(x_{*})$ one needs some continuity-like restrictions for the mapping $F$. This  implication is definitely true for any Hausdorff contraction $F$ and for any upper semicontinuous closed-valued mapping $F$. Moreover, it suffices to assume that the distance function $d_{F}$
is a lower semicontinuous numerical function. So in order to explain the role of the stable positivity  of the distance function $d_{F}$ we shall prove the folowing simple lemma.

\proclaim{Lemma 1.3}
Let $F$ be a closed-valued mapping such that the distance function $d_{F}$ is stably positive. Let $\{x_{n}\}$ be a sequence of points which converges to $x_{*}$, where $x_{n+1} \in F(x_{n})$, for all $n \in \N$. Then $x_{*}$ is a fixed point of $F$.
\endproclaim
\demo {\bf Proof}
Assume to the contrary, i.e. suppose that $x_{*}$ does not belong to the set $F(x_{*})$. Hence $d_{F}(x_{*}) > 0$  and $\inf\{d_{F}(x) \vv  x \in \OO(x_{*})\} = m > 0$ for some neighborhood $\OO(x_{*})$ of the point $x_{*}$. For some number $N \in \N$ onwards, all $x_{n}$ with $n \geq N$ lie in $\OO(x_{*})$. So\, $0 < m \leq d_{F}(x_{n}) \leq d(x_{n}, x_{n+1}) \rightarrow 0, \,\, n \rightarrow \i$. Contradiction.
\q \enddemo

The following theorem is the main result of the paper.

\proclaim{Theorem 1.4} Let numerical functions
$\a:[0,+\i)\rightarrow [1,+\i)$ and $\b:[0,+\i)\rightarrow [0,1)$ be such that:
\itemitem{(1)} $\b(\cdot)$ has the property $(MT)$;
\itemitem{(2)} $\a(t) \leq 1+\g(1-\b(t)),\,\,t\geq 0$, where $\g:(0,1]\rightarrow [0,+\i)$
has the following properties:
\itemitem{(i)} $\g(\cdot)$ is bounded,  $\lim_{s\rightarrow 0} \g(s)=0$;\,\,and
\itemitem{(ii)} the function $p(s)=s-(1-s)\g(s),\,\, 0 \leq s <1$ is essentially positive.

Then every $(\a, \b)-$mapping $F$ of a complete metric space $(X,d)$ into itself is an 
$(\a, \b)-$contraction and has a fixed point, whenever the distance function $d_{F}$ is stably positive.
\endproclaim

In particular, putting $\g(s)=s,\,\,\a(s)=2-\b(s),\,\,p(s)=s^{2}$ in assumptions $(1), (2)$ and equality instead of inequality in assumption $(2)$ we obtain a recent result of \v Ciri\'{c} \cite{2; Theorem 5}. Sometimes there is a possibility for a nonexplicit form of a majorant of the function
$\a:[0,+\i)\rightarrow [1,+\i)$.

\proclaim{Theorem 1.5} Let numerical functions
$\a:[0,+\i)\rightarrow [1,+\i)$ and $\b:[0,+\i)\rightarrow [0,1)$ be such that:
\itemitem{(1)} $\a(\cdot)\b(\cdot)$ has the property $(MT)$; and
\itemitem{(2)} $\a(\cdot)$ is nonincreasing.

Then every $(\a, \b)-$contraction $F$ of a complete metric space $(X,d)$ into itself  has a fixed point, whenever the distance function $d_{F}$ is stably positive.
\endproclaim

The referee has noted that  Theorem 1.5 follows also by \cite{2; Theorem 6}. 
One of the key steps in our proofs of Theorems 1.4 and 1.5 is the "boundary"  property asserting that
$\lim\sup_{s \to 0} \a(s)\b(s) <1$. As a result, it allows us to use
a majorization by a convergent geometrical series. However, sometimes it is also possible to work with the equality $\lim\sup_{s \to 0} \a(s)\b(s) =1$ and hence to use the  property $(R)$ instead of $(MT)$. In the following theorem we present a version with a good power-rate upper estimate, say $\phi$ for the product $\a\b$. In its proof a majorization is made by generalized harmonic series.

\proclaim{Theorem 1.6} Let numerical functions
$\a:[0,+\i)\rightarrow [1,+\i)$ and $\b:[0,+\i)\rightarrow [0,1)$ be such that:
\itemitem{(1)} $\a(\cdot)$ is bounded;\,\,and
\itemitem{(2)} $\a(t)\b(t) \leq \phi(t) =1-Ct^{p}$ for some $C>0,\,\,0<p<1$ and for all $t$ in
some neighborhood of zero.

Then every $(\a, \b)-$contraction $F$ of a complete metric space $(X,d)$ into itself
has a fixed point, whenever the distance function $d_{F}$ is stably positive.
\endproclaim

Finally, in comparison with Theorem 1.4. and  \v Ciri\'{c}'s theorem \cite{2} we have the following:

\proclaim{Example 1.7} There is a finite-valued mapping of the segment $[0,1]$ into itself which:
\itemitem{(1)} is not  a $(2-\b, \b)-$contraction for any $\b:[0,+\i)\rightarrow [0,1)$;
\itemitem{(2)} is not an $(a, \b)-$contraction for any constant $a > 1$; and
\itemitem{(3)} is an $(\a,\b)-$ contraction satisfying all assumptions of Theorem 1.4.
\endproclaim

\head {\bf 2. Proofs }\endhead

We shall organize the proof of Theorem 1.4 in a sequence of  Lemmas 2.1-2.6. The key ingredients are in Lemma 2.5. Its proof preserves the outline of the proof of \cite{2; Theorem 5}.

\proclaim{Lemma 2.1} The product  $\a\b$ has the property $(MT)$.
\endproclaim
\demo {\bf Proof}
Pick any $t \geq 0$. Due to the property $(MT)$ for $\b$ there are numbers $\sigma > 0$ and $0\leq q < 1$ such that
$
\b(s) \leq q, \,\,\, t<s<t+\sigma.
$
Therefore for all such $s$ we have that $1-\b(s)\geq 1-q > 0$ and due to the essential positivity of the function $p(u)=u-(1-u)\g(u)$ we see that $\inf\{p(u)\,|\,\,u\geq 1-q\} = p_{q} > 0$. A  simple calculation
$$
\a(s)\b(s)\leq (1+\g(1-\b(s)))\cdot \b(s) =1-[(1-\b(s))-\b(s)\cdot \g(1-\b(s))]=1-p(1-\b(s)) ,
$$
shows that
$
\a(s)\b(s) \leq 1-p_{q}<1 ,\,\,\,  t<s<t+\sigma.
$
Therefore $\lim\sup_{s \to t+0} \a(s)\b(s) \,\,\leq\,\,1-p_{q} < 1$.
\q \enddemo

\proclaim{Lemma 2.2} If $\e>0$ then $\sup\{\a(t)\b(t) \vv\,\, t \in \a^{-1}([1+\e, +\i))\} = Q_{\e} < 1$.
\endproclaim
\demo {\bf Proof}
Using the equality $\lim_{s\rightarrow 0} \g(s)=0$ we can pick some $\sigma > 0$ such that $\g(s) < \e$ for all $0<s<\sigma$. Then
$$
t \in \a^{-1}([1+\e, +\i))\Longleftrightarrow \a(t)\geq 1+\e  \Longrightarrow 1+\g(1-
\b(t)) \geq 1+\e \Longrightarrow 1-\b(t) \geq \sigma.
$$
Now, exploiting once again the essential positivity of the function $p(\cdot)$, we obtain $\a(t)\b(t) \leq 1-p(1-\b(t)) \leq 1 - \inf\{p(u) \vv\,\,u\geq\sigma\} = Q_{\e} < 1$.
\q \enddemo

\proclaim{Lemma 2.3} For an arbitrary initial point $x_{0} \in X$ there exists a sequence
$\{x_{n}\}_{n=0}^{\i}$ such that for all  $n\in \N$ the point $x_{n+1}$ lies in $F(x_{n})$ and the following properties
hold:
\itemitem{$(A_{n})$} $d_{n}=d(x_{n},x_{n+1}) \leq \a(d_{n})\cdot d_{F}(x_{n})$; and
\itemitem{$(B_{n})$} $d_{F}(x_{n+1}) \leq \b(d_{n})\cdot d_{n}.$
\endproclaim
\demo {\bf Proof}
A straightforward induction using Definition 1.1.
\q \enddemo

\proclaim{Lemma 2.4} Let $\{x_{n}\}_{n=0}^{\i}$ be a sequence constructed in Lemma 2.3.
Then $\{d_{F}(x_{n})\}_{n=0}^{\i}$ is a decreasing numerical sequence, and hence has a nonnegative limit.
\endproclaim
\demo {\bf Proof} Applying Lemmas 2.1 and 2.3 we see that
$$d_{F}(x_{n+1}) \leq \b(d_{n})\cdot d_{n} \leq \b(d_{n})\cdot \a(d_{n})\cdot d_{F}(x_{n})\leq
 (1-p(1-\b(d_{n}))) \cdot d_{F}(x_{n}) < d_{F}(x_{n}). \q $$
 \enddemo

\proclaim{Lemma 2.5} Let $\{x_{n}\}_{n=0}^{\i}$ be a sequence constructed in Lemma 2.3.
Then $\{x_{n}\}_{n=0}^{\i}$ is a fundamental sequence in $X$.
\endproclaim
\demo {\bf Proof}
Clearly, $$d_{F}(x_{n}) = dist(x_{n}, F(x_{n}))\leq d(x_{n},x_{n+1}) = d_{n} \leq C\cdot d_{F}(x_{n}),$$
where $C=\sup \{\a(t) \vv \,\,t\geq 0\} < +\i$.

Denote $\Delta = \lim_{n \to \i} d_{F}(x_{n})$ and $\nabla = \liminf_{n \to \i} d_{n}$. Then
$\Delta \leq \nabla < +\i$. There are exactly three possibilities:
\itemitem{(I)}  $0 < \Delta < \nabla$; \, or
\itemitem{(II)}  $0 < \Delta = \nabla$; \, or
\itemitem{(III)}  $0 = \Delta = \nabla$.

In the case $(I)$ we divide the segment $[\Delta, \nabla]$ into three equal parts. For all sufficiently large indices $n$  the following inequalities hold:
$$
0 < \Delta < d_{F}(x_{n}) < \f{2\Delta + \nabla}{3} < \f{\Delta + 2\nabla}{3} <  d_{n}.
$$

The property $(A_{n})$ (cf. Lemma 2.3) implies that
$$
\a(d_{n}) \geq \f{d_{n}}{d_{F}(x_{n})} > \f{\Delta + 2\nabla}{2\Delta + \nabla} = 1+\f{\nabla - \Delta}{2\Delta + \nabla} = 1 + \e.
$$

By Lemma 2.2 this means that for all sufficiently large indices $n$  the inequality
$\a(d_{n})\b(d_{n}) \leq Q_{\e} < 1$ holds. So
$d_{F}(x_{n+1}) \leq \a(d_{n})\cdot \b(d_{n})\cdot d_{F}(x_{n}) < Q_{\e} \cdot d_{F}(x_{n}).$
Hence from some index, say $N$, onwards the numerical sequence $\{d_{F}(x_{n})\}$  is majorized by a geometrical sequence with the coefficient  $Q_{\e} < 1$. Therefore  $\Delta = 0$. Contradiction.

In the case $(II)$ we have $0 < \nabla = \Delta < d_{F}(x_{n}) \leq d_{n} $ and this is why
$d_{n_{k}} \to \nabla +0,\,\,k \to \i$ for some subsequence. Then $d_{F}(x_{n_{k}}) \to \nabla +0,\,\,k \to \i$, too. By Lemma 2.1 applied to the product $\a\b$ at the point $\nabla = \Delta$, the right upper limit of this product is less than 1. So there exists a number $0 \leq q <1$ such that for all sufficiently large indices $k$ the following inequalities hold:
$$
d_{F}(x_{n_{k+1}}) \leq d_{F}(x_{n_{k}+1}) \leq \a(d_{n_{k}})\cdot \b(d_{n_{k}})\cdot d_{F}(x_{n_{k}}) \leq q\cdot d_{F}(x_{n_{k}}).
$$
Hence from some index, say $K$, onwards the numerical sequence $\{d_{F}(x_{n_{k}})\}$  is majorized by a geometrical sequence with the coefficient with the coefficient $q < 1$. Therefore  $\nabla = 0$. Contradiction.

In the last case $(III)$, recalling that $d_{F}(x_{n}) \leq d_{n} \leq C\cdot d_{F}(x_{n})$, we conclude that $0 = \Delta = \nabla = \liminf_{n \to \i} d_{n} =\lim_{n \to \i} d_{n}$. As in the case $(II)$ one can apply Lemma 2.1 to the product $\a\b$ at the point $0 = \Delta = \nabla$ .
So starting from some index $N$, we have
$$
\sum_{n=N}^{\i} d(x_{n}, x_{n+1})  =  \sum_{n=N}^{\i} d_{n} \leq C\cdot \sum_{n=N}^{\i} d_{F}(x_{n}) \leq C\cdot \sum_{n=N}^{\i} q^{n} < + \i
$$
for some $0 \leq q < 1$. By the triangle inequality, the sequence $\{x_{n}\}$ is fundamental.
\q \enddemo

\proclaim{Lemma 2.6. (End of proof of Theorem 1.4)} Let $\{x_{n}\}_{n=0}^{\i}$ be a sequence, constructed in Lemma 2.3.
Then $\{x_{n}\}_{n=0}^{\i}$ converges to a fixed point of the mapping $F$.
\endproclaim
\demo {\bf Proof}
Completeness of $(X,d)$ together with Lemma 2.5 guarantees the convergence of $\{x_{n}\}_{n=0}^{\i}$ to some point, say $x_{*}$. The inclusion $x_{*} \in F(x_{*})$ has already been checked in Lemma 1.3.
\q \enddemo

\medskip \medskip
\demo {\bf Proof of Theorem 1.5}
As above (cf. Lemma 2.3) for an arbitrary initial point $x_{0} \in X$, there exists a sequence
$\{x_{n}\}_{n=0}^{\i}$ such that for all  $n\in \N$, the point $x_{n+1}$ lies in $F(x_{n})$ and:
\itemitem{$(A_{n})$} $d_{n}=d(x_{n},x_{n+1}) \leq \a(d_{n})\cdot d_{F}(x_{n})$; and
\itemitem{$(B_{n})$} $d_{F}(x_{n+1}) \leq \b(d_{n})\cdot d_{n}.$

Also,
$d_{F}(x_{n+1}) \leq \a(d_{n})\cdot \b(d_{n})\cdot  d_{F}(x_{n}) < d_{F}(x_{n}), $
whereas
$d_{n+1} \leq \a(d_{n+1})\cdot d_{F}(x_{n+1}) \leq  \a(d_{n+1})\cdot
\b(d_{n})\cdot d_{n}. $
Hence, as above (cf. Lemma 2.4) $\{d_{F}(x_{n})\}_{n=0}^{\i}$ is a decreasing sequence which converges to some $\Delta \geq 0$. It turns out that $\{d_{n})\}_{n=0}^{\i}$ is a decreasing sequence, too.
Indeed, assume to the contrary that $d_{n+1} \geq d_{n}$ for some $n \in \N$. \, Then $\a(d_{n+1}) \leq \a(d_{n})$ and
$$
d_{n} \leq d_{n+1} \leq  \a(d_{n+1})\cdot \b(d_{n})\cdot d_{n} \leq \a(d_{n})\cdot \b(d_{n})\cdot d_{n} < d_{n}.
$$
Contradiction. So $d_{n} \rightarrow \nabla+0,\,\,n \to \i$.

The inequality $d_{F}(x_{n+1}) \leq \a(d_{n})\cdot \b(d_{n})\cdot  d_{F}(x_{n})$,\, together with the assumption $(1)$ and  $d_{n} \rightarrow \nabla+0,\,\,n \to \i$,\, immediately implies that $\Delta=0$. It now follows from the inequality $d_{F}(x_{n}) \leq d_{n} \leq \a(0) \cdot d_{F}(x_{n})$  that $\nabla=0$. The rest of the proof is the same as the proof of Theorem 1.4.
\q \enddemo

\medskip
\demo {\bf Proof of Theorem 1.6}
We preserve the choice of $x_{n+1} \in F(x_{n})$ as in the proofs above (cf. Lemma 2.3). So
$$
d_{F}(x_{n+1}) \leq \a(d_{n})\cdot \b(d_{n})\cdot d_{F}(x_{n})\leq \phi(d_{n})\cdot d_{F}(x_{n})
\leq \phi(d_{F}(x_{n}))\cdot d_{F}(x_{n}),
$$
because $d_{F}(x_{n}) \leq d_{n}$ and the majorant $\phi(\cdot)$ is a decreasing function.

Henceforth, in order to check the convergence of the series $\sum_{n=0}^{\i} d_{F}(x_{n})$
it suffices to show that the series $\sum_{n=0}^{\i} \phi_{n}(t)$ is convergent for all $t>0$, where
$$
\phi_{0}(t)=t, \qq \phi_{n+1}(t)= \phi(\phi_{n}(t)) \cdot \phi_{n}(t).
$$

In other words, by fixing $t$ and omitting $t$ in the brackets, we only need to show
that $\sum_{n=0}^{\i} \phi_{n} < +\i$. We shall complete the proof by checking that
$$
\phi_{n} \leq \f{1}{(C_{1} n + C_{2})^{\f{1}{p}}} \qq \Longleftrightarrow \qq
\f{1}{\phi_{n}^{p}} \geq C_{1} n + C_{2}
$$
for some constants $C_{1}>0, C_{2}>0$ and for all sufficiently large $n$.
Let us verify that one can take $C_{1}=pC,\,\, C_{2}=\phi_{0}^{-p}$. To this end we represent
$$
\f{1}{\phi_{n}^{p}} = \left( \f{1}{\phi_{n}^{p}} - \f{1}{\phi_{n-1}^{p}} \right)
+...+\left(\f{1}{\phi_{1}^{p}} - \f{1}{\phi_{0}^{p}}\right) + \f{1}{\phi_{0}^{p}}.
$$

Next, for each $k=0,1,...,n-1$, by using the Lagrange theorem we see that
$$
\f{1}{\phi_{k+1}^{p}} - \f{1}{\phi_{k}^{p}}=\f{1}{(\phi(\phi_{k}) \cdot
\phi_{k})^{p}} - \f{1}{\phi_{k}^{p}} = \f{1}{\phi_{k}^{p}} \cdot \f{1-(1-C\phi_{k}^{p})^{p}}{(1-C\phi_{k}^{p})^{p}}=\left( s=C \cdot \phi_{k}^{p} \right)=
$$

$$
=C \cdot \f{1}{s} \cdot \f{(-s)\cdot p \cdot (1-\tau)^{p-1} \cdot (-1)}{(1-s)^{p}}=
pC \cdot \left( \f{1-\tau}{1-s} \right)^{p} \cdot \f{1}{1-\tau} > pC=C_{1}
$$
because $0<\tau < s$.

Hence the series $\sum_{n=0}^{\i} d_{F}(x_{n})$ is convergent and
the series $\sum_{n=0}^{\i} d_{n}$ is also convergent, because
$d_{F}(x_{n}) \leq d_{n} \leq \sup \{\a(t): t \geq 0 \} \cdot d_{F}(x_{n}) $. Finally, we see that
the sequence $\{x_{n}\}$ is fundamental. The rest of the proof is standard.
\q \enddemo

\medskip

Now we pass to Example 1.7 and we use the idea of \cite{2,7} except that our construction will avoid rather unexpected constants as $\f{15}{32}, \f{7}{24}, \f{17}{96}$, etc., and will be based only on piecewise linear functions. We define an almost singlevalued mapping $F:[0,1] \to [0,1]$ by setting
$F(1)=\{\f{1}{3},\f{3}{4}\}$ and $F(x)=\{\f{2x}{3}\}$ if $0 \leq x < \f{3}{4}$,\,$F(x)=\{\f{1}{2}\}$\,otherwise.

\demo {\bf Proof of (1)}
Suppose to the contrary that the mapping $F$ is an $(2-\b, \b)-$contraction for some numerical function $\b:[0,+\i)\rightarrow [0,1)$. Consider the point $x=1$.

If one chooses the point $y=\f{3}{4} \in F(1)=F(x)$ then $d(x,y)=\f{1}{4}$ and for the condition $(B)$ from Definition 1.1 one obtains
$d_{F}(y)=dist(y,F(y))=d\left(\f{3}{4}, \f{1}{2}\right)=\f{1}{4}\leq \b\left(\f{1}{4}\right)\cdot \f{1}{4},
$
or, $\b(\f{1}{4}) \geq 1$. Contradiction.

For the other choice $y=\f{1}{3} \in F(1)=F(x)$ and for the condition $(A)$ from Definition 1.1. we conclude that
$
d(x,y)=\f{2}{3} \leq \a\left(\f{2}{3}\right)\cdot d_{F}(x)=\a\left(\f{2}{3}\right)\cdot \f{1}{4},
$
or, $\a(\f{2}{3})\geq \f{8}{3}>2$, which contradicts the fact that $\a=2-\b < 2$.
\q \enddemo

\demo {\bf Proof of (2)}
Suppose to the contrary that the mapping $F$ is an $(a, \b)-$contraction for some constant $a > 1$ and for some numerical function
$\b:[0,+\i)\rightarrow [0,1)$. Consider the point $x=\f{1}{2}$. Then $F(x)= \{\f{1}{3}\}$ and there is a unique choice of $y \in F(x)$, namely $y=\f{1}{3}$. Hence,
$F(y)= \{\f{2}{9}\}$ and $d(x,y)=d_{F}(x)=\f{1}{6}$ and $d_{F}(y)= \f{1}{3} - \f{2}{9}=\f{1}{9}$. So the condition $(B)$ looks as follows:
$
d_{F}(y) \leq \b(d(x,y))\cdot d(x,y) \Longleftrightarrow \b\left(\f{1}{6}\right)\geq \f{2}{3}.
$

As in the previous proof, for the point $x=1$ we have $\b(\f{1}{4}) \geq 1$ for the chosen point $y=\f{3}{4}$, \, or we have $a=\a(\f{2}{3})\geq \f{8}{3}$ for the chosen point $y=\f{1}{3}$. But $a  < \f{3}{2}$ because of the restriction $\a(\f{1}{6})\cdot \b(\f{1}{6}) = a\cdot \b(\f{1}{6}) < 1$.\,Contradiction.
\q \enddemo

\demo {\bf Proof of (3)}
Define $\a(\cdot)$ and $\b(\cdot)$ by the setting $\a(x)=\f{4}{3}$ for $0\leq x \leq \f{1}{2}$,\,\, $\a(x)= \f{8}{3}$ for $\f{1}{2} < x$ and $\b(x)=\f{2}{3}$ for $0\leq x \leq \f{1}{3}$,\,\, $\b(x)=\f{1}{2}$ for $\f{1}{3} < x \leq \f{1}{2}$, \,\,$\b(x)=\f{1}{3}$ for $\f{1}{2} < x$.

Then the function $\g(\cdot)$ with the property that $\a(x)=1+\g(1-\b(x))$ can be defined by $\g(\f{1}{3})=
\f{1}{3}=\g(\f{1}{2}),\,\,\g(\f{2}{3})=\f{5}{3}\,\,$ and $\g(t)=0$ otherwise.
Therefore for $p(x)=x-(1-x)\g(x)$ we have $p(\f{1}{3})=
\f{1}{9}=p(\f{2}{3}),\,\,p(\f{1}{2})=\f{1}{3},\,\,p(x)=x$ otherwise and thus all assumptions $(1),(2)$ of Theorem 1.4  are satisfied.

Let us check that $F$  is really an $(\a, \b)-$contraction. For the point $x=1$ we choose $y= \f{1}{3} \in F(x)$. Then $F(y)=\left\{\f{2}{9} \right\}$. The inequality $d(x,y)\leq \a(d(x,y))\cdot d_{F}(x)$ from the condition $(A)$ becomes $\f{2}{3} \leq \a(\f{2}{3}) \cdot \f{1}{4}$ and is true because $\a(\f{2}{3})= \f{8}{3}$. Also,  the inequality $d_{F}(y)\leq \b(d(x,y))\cdot d(x,y)(x)$ from the condition $(B)$ becomes $\f{2}{9} \leq \b(\f{2}{3}) \cdot \f{2}{3}$ and it holds because $\b(\f{2}{3})= \f{1}{3}$.

For the point $0 \leq x <1$ there are unique $y \in F(x), \,\,y=\f{2x}{3}$ and $z \in F(y), \,\,z=\f{4x}{9}$. In particular, $d(x,y)=d_{F}(x)$ and the inequality $d(x,y)\leq \a(d(x,y))\cdot d_{F}(x)$ from the condition $(A)$ holds because $\a(\cdot) > 1$. For checking of $(B)$ we calculate
\medskip
$d(x,y)=x-\f{2x}{3}=\f{x}{3}, \,\,0\leq x \leq \f{3}{4},\qq d(x,y)=x-\f{1}{2}, \,\,
\f{3}{4}\leq x \leq 1,
$
\medskip
$d_{F}(y)=\f{2x}{3}-\f{4x}{9}=\f{2x}{9}, \,\,0\leq x \leq \f{3}{4},\qq d_{F}(y)=\f{1}{6}, \,\,
\f{3}{4}\leq x \leq 1.
$
\medskip
So for $0\leq x \leq \f{3}{4}$, the inequality $d_{F}(y)\leq \b(d(x,y))\cdot d(x,y)$ from $(B)$ becomes $\f{2x}{9} \leq \b(\f{x}{3}) \cdot \f{x}{3} = \f{2}{3} \cdot \f{x}{3}$.
For points $\f{3}{4} < x \leq \f{5}{6}$ we have
$\f{1}{4} < d(x,y)=x-\f{1}{2} \leq \f{1}{3}$ and the needed inequality  becomes
$\f{1}{6} = \f{2}{3} \cdot \f{1}{4} < \f{2}{3} \cdot (x-\f{1}{2})=\b(x-\f{1}{2}) \cdot (x-\f{1}{2}) $. Finally, if $\f{5}{6} < x <1$ then $\f{1}{3} < d(x,y)=x-\f{1}{2} < \f{1}{2}$ and and the desired inequality becomes
$\f{1}{6} = \f{1}{2} \cdot \f{1}{3} < \f{1}{2} \cdot (x-\f{1}{2}) = \b(x-\f{1}{2}) \cdot (x-\f{1}{2}) $.
\q \enddemo

\head {\bf 3. Concluding remarks}\endhead

\noindent
{\bf Remark 3.1.} Lower semicontinuity of a nonnegative numerical function implies its stable positivity, but not vice versa. For example, let $F$ be a Hausdorff $k-$contraction of $(X,d)$ and  $x_{0}$  a non-fixed point of $F$. Define a new mapping, say $G$, by setting $G(x)=F(x),\,\,x \not= x_{0}$ and letting $G(x_{0})$ be an arbitrary closed subset of $X$ with $dist(x_{0}, F(x_{0})) < dist(x_{0}, G(x_{0}))$. Then the distance function $d_{G}(\cdot)$ is evidently stably positive over $X$, but $d_{G}(\cdot)$ is not lower semicontinuous at the point $x_{0}$. Observe also that the upper semicontinuity of a closed-valued mapping $F$ implies the lower semicontinuity of the distance function $d_{F}$ which, in turn, implies its stable positivity, due to Lemma 1.3.

\noindent
{\bf Remark 3.2.} As it was pointed above, the special case of Theorem 1.4 with $\g(s)=s,\,\,\a(s)=2-\b(s),\,\,p(s)=s^{2}$ in fact coincides with a recent result of \v Ciri\'{c}
\cite{2; Theorem 5}. 
For the mappings with proximinal values it suffices to take $\g(s)=0,\,\,\a(s)=1,\,\,p(s)=s$ to obtain results of \cite{4,7} and \cite{2; Theorem 7}. In particular, we obtain results on fixed points for singlevalued contractions and for compact-valued Hausdorff contractions. Mizoguchi-Takahashi theorem \cite{8} is also a special case of Theorem 1.4, as was discussed above
(after Definition 1.1).

\noindent
{\bf Remark 3.3.} The essential positivity of a numerical function $h:[0,\i) \to [0;1)$ is equivalent to the property that $h$ admits an increasing positive minorant. Hence the assumption $(2.ii)$ of Theorem 1.4, after dividing by $1-s$, can be expressed in the form
$$
\g(s)=\f{s}{1-s} - \mu(s) \qq \Longleftrightarrow \qq \g(s)=(s+s^{2}+s^{3}+...)-\mu(s)
$$
for some positive increasing function $\mu(\cdot)$.

In particular, if one  takes an arbitrary function $\g:(0,1]\rightarrow [0,+\i)$ such that $0 < \g(s)\leq s+s^{2}+...+s^{m},\,\,m \in \N$ then the assumption $(2.i)$ from
Theorem 1.4 is evident, whereas $(2.ii)$ holds because
$$
p(s)=s-(1-s)\g(s) \geq s-(1-s)(s+s^{2}+...+s^{m})=s^{m+1}.
$$
Hence, these cases  are the sources of new fixed point theorems.

\noindent
{\bf Remark 3.4.} Recall that Reich \cite{9} showed that the property $(R)$ guarantees  existence of fixed points for an arbitrary compact-valued $k(\cdot)$-contraction $F$
of a complete metric space. In 1974 he also proposed  \cite{10} the still unresolved
problem on  possibility of removing the compactness condition. In 1989 Mizoguchi and Takahashi \cite{8} obtained the result for any closed-valued $k(\cdot)$-contractions but under the condition $(MT)$ stronger than $(R)$. So we still have the following interesting open problem.

\proclaim{Question 3.5} Is Theorem 1.4 true for compact-valued $(\a,\b)-$contractions but with $(MT)$ replaced by $(R)$ in the assumption $(1)$?
\endproclaim

Roughly speaking, the difficulty is that Reich proved his theorem by passing to the compact exponent of $(X,d)$ endowed with the Hausdorff distance $H_{d}$  and by using an appropriate fixed point theorem for singlevalued mapping of such complete metric space into itself. 
However, in our case a compact-valued $(\a,\b)-$contraction does not generate a singlevalued mapping of the compact exponent into itself. Moreover, in the absence of Hausdorff distance  one needs to find another kind of "metric" in the compact exponent which agrees with the notion of an $(\a,\b)-$contraction.

\noindent
{\bf Remark 3.6.} Formally, Theorem 1.6 admits the following abstract form.
\proclaim{Theorem 3.7} Let numerical functions
$\a:[0,+\i)\rightarrow [1,+\i)$ and $\b:[0,+\i)\rightarrow [0,1)$ be such that:
\itemitem{(1)} $\a(\cdot)$ is bounded;\,\,and
\itemitem{(2)} the product $\a(t)\b(t)$ has a nonincreasing
majorant $\phi(t) < 1$ with  $\sum_{n=0}^{\i} \phi_{n}(t) <+\i$, where
$
\phi_{0}(t)=t,  \phi_{n+1}(t)= \phi(\phi_{n}(t)) \cdot \phi_{n}(t).
$

Then every $(\a, \b)-$contraction $F$ of a complete metric space $(X,d)$ into itself
has a fixed point, whenever the distance function $d_{F}$ is stably positive.
\endproclaim

A special case of Theorem 1.6 for Hausdorff $k(\cdot)-$contractions with power-rate majorants for $k(\cdot)$ yields the main result of \cite{3}. See also \cite{11} for examples of (in this sense) "summable" functions $\phi(\cdot)$ which have no power-rate upper estimates.

\head {\bf Acknowledgements}\endhead

This research was supported by Slovenian Research Agency grants P1-0292-0101 and J1-2057-0101. The author thanks Pavel V. Semenov and the referee for comments and suggestions.

\enddocument